\documentclass[11 pt, psamsfonts]{amsart}

\usepackage{amssymb}

\calclayout

\def\mathbold{\pmb}

\def\ignore#1{\relax}

\def\th{^{{\rm th}}}
\def\ncps{algebraic noncommutative probability space}
\def\ngraded{$\Z_n$--graded }
\def\nc{{\rm NC}}

\newcommand\newsect{\setcounter{equation}{0}\section}

\newcommand\be{\begin{equation}}
\newcommand\ee{\end{equation}}

\def\Z{{\mathbb Z}}

\def\C{{\mathbb C}}

\def\A{\mathcal A}
\def\B{\mathcal B}
\def\L{\mathcal L}
\def\P{\mathcal P}

\def\one{\mathbf 1}

\def\inv{^{-1}}

\def\la{\lambda}

\def\cardinality #1 {{\rm card}(#1)}
\def\card{\cardinality}
\def\th{^{{\rm{th}}}}
\def\r #1 #2{R^{(#1)}_{#2}}
\def\a #1 #2{\alpha_{#2}^{(#1)}}
\def\rest #1{_{| #1}}
\def\Ad{{\rm Ad}}

\theoremstyle{plain}
\newtheorem{theorem}{Theorem}[section]

\theoremstyle{plain}
\newtheorem{proposition}[theorem]{Proposition}

\theoremstyle{plain}
\newtheorem{corollary}[theorem]{Corollary}

\theoremstyle{plain}
\newtheorem{lemma}[theorem]{Lemma}

\theoremstyle{definition}
\newtheorem{definition}[theorem]{Definition}
\theoremstyle{definition}
\newtheorem{example}[theorem]{Example}

\theoremstyle{exercise}

\theoremstyle{definition}
\newtheorem{remark}[theorem]{Remark}

\theoremstyle{remark}

\begin{document}

\title{
$\Z_n$--graded independence }

\author{ \bf Frederick M. Goodman}

\address{ Department of Mathematics\\ University of Iowa\\ Iowa
City, Iowa}

\email{ goodman@math.uiowa.edu}

\thanks{I am grateful to Andu Nica and Florin Radulescu for
suggesting the problem of developing the linearizing transform for
\ngraded independence.}
\date{June 27, 2002}

\maketitle

\newsect{Introduction.}
\medskip

This note extends the work of J. Mingo and A. Nica  ~\cite{Mingo
Nica} on graded independence of random variables and the graded
$R$--transform.  We generalize the results of Mingo and Nica
from the
context of 
$\Z_2$--graded noncommutative probability spaces to that of
$\Z_n$--graded noncommutative probability spaces.

A very general setting for noncommutative probability is an 
{\em algebraic noncommutative probability space} consisting of a not
necessarily commutative unital complex algebra $\A$ and a linear
functional
$\varphi : \A \rightarrow \C$ satisfying $\varphi(\one) = 1$. 
A noncommutative random variable is simply an element $x$ of the algebra
$\A$, and its expectation is $\varphi(x)$.

The general question addressed in ~\cite{Mingo Nica} is the existence of
notions of independence of random variables other than (and possibly
interpolating between) the classical  concept of independence  and the
concept of free independence of Voiculescu ~\cite{VDN}.  
Several authors have provided axiomatizations and classification
theorems for independence on   algebraic noncommutative probability
spaces (~\cite{Speicher Universal}, ~\cite{Muraki}, ~\cite{Muraki2},
~\cite{GS}, ~\cite{Franz what is independence}).   The common conclusion
of these works (under various axiomatizations) is that the only
possibilities are classical
 independence and free independence.  (There are additional
possibilities, namely, monotone and boolean independence, for
nonunital algebras.)

Nevertheless, for more restricted categories of algebras, with additional
structure, there can be notions of independence specific to these
categories, which take into account the additional structure.
Here we consider {\em $\Z_n$--graded independence} of {\em
$\Z_n$--graded  noncommutative probability spaces.}

\subsection{\ngraded independence}
\medskip
\begin{definition}  An  $\Z_n$--{\em graded algebraic noncommutative
probability space}
$(\A, \varphi, \gamma)$ consists of
\begin{enumerate}
\item
a unital complex algebra $\A$,  
\item
a linear functional $\varphi$
on $\A$ satisfying $\varphi(\one) = 1$, and
\item
a unital algebra automorphism $\gamma$ of $\A$ of order $n$ satisfying
$\varphi
\circ \gamma = \varphi$.
\end{enumerate}
A {\em morphism} $\Phi : (\mathcal A_1, \varphi_1, \gamma_1) \rightarrow
(\mathcal A_2, \varphi_2, \gamma_2) $ of $\Z_n$--graded
algebraic noncommutative probability spaces is a unital algebra morphism
satisfying
$\gamma_2\circ \Phi = \Phi\circ\gamma_1$  and
$\varphi_2 \circ \Phi = \varphi_1$.
\end{definition}
\medskip

If $(\A, \varphi, \gamma)$ is a 
$\Z_n$--graded algebraic noncommutative
probability space, then $\A$ is the direct sum of
eigenspaces for $\gamma$,
$$
\A = \bigoplus_{q \in \Z_n^*} \A_q,
$$
where $ \Z_n^*  $  denotes the group of $n\th$ roots of unity in
$\C$, and  for $q \in  \Z_n^*  $,  
$$\A_q = \{ a \in \A : \gamma(a) = qa\}.
$$
 The elements of any eigenspace
are said to be {\em homogeneous}.  The eigenspaces satisfy 
$$
A_q A_r \subseteq A_{q r} \quad \text{for }  q, r \in \Z_n^*.
$$ 
The  eigenspace
$A_1$ is the {\em subalgebra of fixed points for $\gamma$}.  Moreover, for
$r
\ne 1$, the restriction of $\varphi$ to $\A_r$ is identically zero,
since for $x \in A_r$, we have $\varphi(x) = \varphi(\gamma(x)) =
r \varphi(x)$.  Equivalently, $\varphi = \varphi \circ E$,
where $E$ is the projection of $\A$ onto $\A_1$ defined
by 
$$E(x) = (1/n) \displaystyle\sum_{0 \le i \le n-1} \gamma^i(x).
$$

Let $q$ be any primitive $n\th$ root of unity, so $\Z_n^*$ is generated
by $q$ as a multiplicative group.  If
$a \in
\A_{q^r}$, we say the degree of
$a$  with respect to $q$ is
$r$, and we write
$\delta_q(a) = r$.  For homogeneous elements $a$ and $b$, we have
$\delta_q(ab) \equiv \delta_q(a) + \delta_q(b)  \mod n$.

\medskip
\begin{definition}
 Let $\A$ and $\mathcal B$
be two unital subalgebras of a $\Z_n$--graded noncommutative probability
space $(\mathcal C, \varphi, \gamma)$, with $\A$ and $\mathcal B$   
invariant under the grading automorphism
$\gamma$. The algebras $\A$ and $\mathcal B$ are said to be
{\em$\Z_n$--graded independent} if 
\begin{enumerate}
\item $\A$ and $\mathcal B$ gradedly commute; that is, there
exists a primitive $n\th$ root of unity $q$ such that if
for all $a \in \A$ and $b \in \mathcal B$,   
$$
b a =
q^{\delta_q(a)\delta_q(b)} a b.
$$
\item  $\varphi(a b) = \varphi(a)\varphi(b)$ for all
$a \in \A$ and $b \in \mathcal B$. 
\end{enumerate} 
\end{definition}

\medskip
Two examples of algebras with $\Z_n$--graded independent subalgebras are
the  {\em $\Z_n$--graded tensor product} of two 
$\Z_n$--graded graded algebras and the {\em rotation  algebras} with
parameter an $n\th$ root of unity.  See Section \ref{section examples} 
for a discussion of these examples.

\medskip
\subsection{Linearization}
In both classical and free probability, there exist transforms which
linearize the addition of independent random variables.  
Let $(\A, \varphi)$ be a noncommutative probability space.
For $a \in \A$, the {\em moment sequence of} $a$ is the sequence
of numbers
$(\mu_k = \varphi(a^k))_{k \ge 0}$.   The moment sequence 
determines a linear
functional
$\mu_a : \C[X] \rightarrow \C$  by  $\mu_a(X^k) = \mu_k$ for $k \ge
0$.   A {\em linearizing transform} is a bijective correspondence
$(\mu_k)_{k \ge 0} \mapsto \alpha(a) = (\alpha_k)_{k \ge 1}$ such that
for independent random variables $a, b$, we have
$\alpha(a + b) = \alpha(a) + \alpha(b)$.  The sequence
$\alpha$ is called the {\em cumulant} sequence.
To obtain explicit formulas and to facilitate computation, it is useful
to use introduce an (exponential or ordinary) generating function
$R(z) = R[\mu_a](z) $ for  the cumulant sequence $\alpha(a)$.  The
linearizing property of the transform is then expressed as
$$
R[\mu_{a + b}](z) = R[\mu_a](z) + R[\mu_b](z),
$$
for independent random variables $a$ and $b$.

\subsubsection{Classical cumulants and classical independence} 
Let $\mu: \C[x] \rightarrow \C$ be a unital linear functional, with
moments $\mu_k = \mu(X^k)$, $k \ge 0$.  The {\em classical cumulant}
sequence corresponding to $\mu$ is determined recursively by the
relations
\begin{equation} \label{equation definition classical cumulants}
\mu(X^k) = \sum_{P \in \P[k]}  \prod_{B \text{ a block
of } P}
\alpha_{\card {B} } \qquad (k \ge 1),
\end{equation}
where $ \P[k]$ denotes the family of set partitions of the interval
$[k] = {1, 2, \dots, k}$.  For example,
\begin{equation*}
\begin{aligned}
\mu_1 &= \alpha_1 \\
\mu_2 &=  \alpha_1^2 + \alpha_2  \\
\mu_3 &=  \alpha_1^3 + 3 \alpha_1 \alpha_2 + \alpha_3\\
\mu_4 &= \alpha_1^4 + 4 \alpha_1 \alpha_3 + 3 \alpha_2^2 + 6 \alpha_2
\alpha_1^2 +\alpha_4.
\end{aligned}
\end{equation*}
The exponential generating function for the classical cumulants
$$
R_1[\mu](z) = \sum_{k \ge 1} \frac{\alpha_k}{k!} z^k
$$
satisfies
\begin{equation} \label{equation formula classical R transform}
R_1[\mu](z) = \log\left[ \sum_{k \ge 0} \frac{\mu_a(X^{k })}
{k !} \  z^{k }  
\right].
\end{equation}
See, for example, Section II.12.8 in ~\cite{Shiryayev}.
For classically independent random variables $x$ and $y$, we have
\begin{equation}
R_1[\mu_{x+y}](z)  = R_1[\mu_x](z) + R_1[\mu_y](z)
\end{equation}
Formally, we have $R_1[\mu_x](z) = \log \varphi \exp(z x)$, and
$$
\begin{aligned}
R_1[\mu_{x+y}](z) &= \log[\varphi \exp(z (x + y))    ]  
= \log[\varphi(\exp(z x) \exp(z y) )]  \cr
&= \log [\varphi(\exp(z x))  \varphi(\exp(zy))] \cr &=
\log[\varphi[\exp(zx)] + \log[\varphi[\exp(zy)] \cr
&= R_1[\mu_{x}](z) +  R_1[\mu_{y}](z).
\end{aligned}
$$
These computations can be justified by formal power series manipulations.

\subsubsection{Free cumulants and free independence}
The {\em free cumulants} of a moment sequence $\mu_k = \mu(X^k)$ are
determined recursively by the requirements
\begin{equation} \label{equation definition free cumulants}
\mu(X^k) = \sum_{P \in \nc[k]}  \prod_{B \text{ a block
of } P}
\alpha_{\card {B} } \qquad (k \ge 1),
\end{equation}
where now the sum is over {\em noncrossing partitions} of 
$[k]$.  Free cumulants were introduced by Speicher ~\cite{Speicher 1}.
For example,
\begin{equation*}
\begin{aligned}
\mu_1 &= \alpha_1 \\
\mu_2 &=  \alpha_1^2 + \alpha_2  \\
\mu_3 &=  \alpha_1^3 + 3 \alpha_1 \alpha_2 + \alpha_3\\
\mu_4 &= \alpha_1^4 + 4 \alpha_1 \alpha_3 + 2 \alpha_2^2 + 6 \alpha_2
\alpha_1^2  +\alpha_4.
\end{aligned}
\end{equation*}
The ordinary generating function for the free cumulants 
\begin{equation*}
R_0[\mu](z) = \sum_{k \ge 1} {\alpha_k} z^k
\end{equation*}
is the $R$--transform of Voiculescu ~\cite{Voi R transform}, in the
formulation due to Speicher ~\cite{Speicher 1}.  The $R$--transform
satisfies
\begin{equation*}
\sum_{k \ge 0}\frac{\mu(X^k)}{\zeta^k} = \left[ \frac{1 +
R_0[\mu](z)}{z}  
\right]^{<-1>} (\zeta),
\end{equation*}
where $<-1>$ indicates inversion with respect to composition of power
series.  $R_0$   linearizes addition of free--independent random
variables,
\begin{equation*}
R_0[\mu_{a + b}](z) = R_0[\mu_{a }](z)  + R_0[\mu_{ b}](z),
\end{equation*}
when $a$ and $b$ are free--independent random variables.  See
~\cite{Voi R transform} and ~\cite{Speicher 1} for details and proofs.

\subsubsection{$q$-cumulants and graded independence} 
Nica and Mingo (~\cite{Nica cmp} and ~\cite{Mingo Nica})
observed that there exists a simultaneous $q$--deformation of
the formulas (Equations (\ref{equation definition classical cumulants})
and  (\ref{equation definition free cumulants}))  defining the 
classical and free cumulants.  For a parameter $q$, they define the
{\em $q$--cumulants} as follows.

\medskip
\begin{definition}
 The {\em $q$-cumulants}  $(\a q k)_{k \ge 1}$ of a moment sequence
$\mu_k = \mu(X^k)$ are determined recursively by the requirements
\begin{equation}
\mu(X^k) = \sum_{P \in \P[k]} q^{c_0(P)} \prod_{B \text{ a block
of } P}
\a q {\card {B} } \qquad (k \ge 1),
\end{equation}
where $c_0(P)$ is the {\em reduced crossing number of the partition $P$},
defined in  Definition
\ref{definition reduced crossing number}.
\end{definition}
\medskip

For example,
\begin{equation*}
\begin{aligned}
\mu_1 &= \a q 1 \\
\mu_2 &=  (\a q 1)^2 + \a q 2  \\
\mu_3 &=  (\a q 1)^3 + 3 \a q 1 \a q 2 + \a q 3\\
\mu_4 &= (\a q 1)^4 + 4 \a q 1 \a q 3 + (2 + q) (\a q 2)^2    + 6
\a q 2(\a q 1)^2 +
\a q 4.
\end{aligned}
\end{equation*}
For $q = 1$, one recovers the classical cumulants, while for 
$q = 0$, one obtains the free cumulants, since $c_0(P) = 0$ precisely when
$P$ is noncrossing.

 Mingo and Nica asked whether $\mu \mapsto \a q {}$
is a linearizing transform for some species of independent random
variables.  They showed that, in fact, for $q = -1$, 
$\mu \mapsto \a {-1} {}$ linearizes addition of $\Z_2$--graded
independent random variables.  {\em The main result of the present paper
is that for $q$ a primitive
$n\th$ root of unity, $\mu \mapsto \a q {}$ linearizes addition of
 $\Z_n$--graded
independent random variables.}

For $q$ a proper, but not necessarily primitive, $n\th$ root of
unity,  we define the $\r n q$--transform of a moment sequence $\mu_k =
\mu(X^k)$ by
$$
\r n q[\mu](z) = \sum_{k \ge 1} \frac{\a q {nk} }{k !} \ z^{nk}.
$$
Moreover, for $a$ a homogeneous element of a \ngraded \ncps, and for
$q$ a {\em primitive} $n\th$ root of unity we define the graded
$\r q n$ transform of $a$ to be
$$
\r n q[a](z) = \r  {n'}  {q^{\delta(a)}}[\mu_a](z),
$$
where $n' = n/{\gcd(\delta_q(a),n)}$ and $\delta(a)$ is the $q$--degree
of $a$.

We show that
$$
\r n q[a](z) = \log\left[ \sum_{k \ge 0} \frac{\mu_a(X^{k n'})}
{k !} \  z^{k n'}  
\right],
$$
and
$$
\r  n q [a_1 + a_2] = \r n q [a_1]  + \r n q [a_2] 
$$
for homogeneous \ngraded independent random variables of the same degree
(Theorems \ref{theorem formula for rq transform}  and 
\ref{theorem rq transform linearizes}).

Our approach closely
follows that of ~\cite{Mingo Nica}, but several technical innovations are
necessary for the proofs.

\newsect{Examples of \ngraded independence} \label{section examples}

\subsection{Graded tensor product} The $\Z_n$--graded tensor product
of two $\Z_n$--graded algebraic noncommutative probability spaces
depends on the choice of a primitive $n\th$ root of unity $q$.
Let 
 $(\mathcal A, \varphi, \alpha)$ and  
$(\mathcal B, \psi, \gamma) $  be two
$\Z_n$--graded algebraic noncommutative probability spaces.
Define $\mathcal A \otimes_q \mathcal B$ to be 
$\mathcal A \otimes \mathcal B$ as a vector space, with the
algebra structure given by
$$
(a \otimes b) (a' \otimes b') = q^{\delta_q(b) \delta_q(a')} a a' \otimes
b b'.
$$
Then $(\A \otimes_q \mathcal B, \varphi \otimes \psi, \alpha
\otimes
\gamma)$  is a $\Z_n$--graded algebraic noncommutative probability space,
and the injections $\A \rightarrow \A \otimes_q \mathcal B$ 
and $\mathcal B \rightarrow \A \otimes_q \mathcal B$  given\
by $a \mapsto a \otimes_q 1$, $b \mapsto 1 \otimes b$ are
morphisms of $\Z_n$--graded algebraic noncommutative probability spaces.

The tensor product $\otimes_q$ is associative, but not commutative.

A more general construction of graded tensor products has been discussed,
for example, in ~\cite{Knus graded}.  The $\Z_n$ graded tensor product has
the following universal property:

\medskip

\begin{proposition} Suppose
\begin{enumerate}
\item  $\A, \mathcal B, \mathcal A',\mathcal B',\mathcal  C$ are
$\Z_n$--graded algebras with $\mathcal A', \mathcal B' \subseteq C$.  
\item  \ $\varphi : \A \rightarrow \mathcal A'$ and
$\psi : \mathcal B \rightarrow \mathcal B'$ are morphisms of
$\Z_n$--graded algebras.
\item $q$ is a primitive $n\th$ root of unity.
\item $\mathcal A'$ and $ \mathcal B'$ \ $q$-commute in $ \mathcal C$;
that is, for homogeneous elements
$a \in \mathcal A'$ and $b \in \mathcal B'$,
$$
b a = q^{\delta_q(a) \delta_q(b) } a b.
$$
Then there is a unique morphism of $\Z_n$--graded algebras
$\varphi \otimes_q \psi : \A \otimes_q \mathcal B \rightarrow \mathcal C$
such that
$(\varphi \otimes_q \psi)(a\otimes 1) = \varphi(a)$ for all $a \in \A$
and $(\varphi \otimes_q \psi)(1\otimes b) = \psi(b)$ for all $b \in
\mathcal B$.
\end{enumerate}
\end{proposition}

\begin{proof}  Left to the reader.
\end{proof}

\medskip

\begin{proposition}  Let $(\mathcal C, \varphi, \gamma)$ be a \ngraded
\ncps, and let $\A, \mathcal B$ be subalgebras of $\mathcal C$ invariant
under
$\gamma$.  Suppose $q$ is a primitive $n\th$ root of unity such that
$\A$ and $\mathcal B$ \ $q$-commute.  The following are equivalent:
\begin{enumerate}
\item  $\A$ and $\mathcal B$  are \ngraded independent.
\item There is a morphism of \ngraded algebraic noncommutative probability spaces
 $$h:(\A \otimes_q \mathcal B, \varphi\otimes \varphi)  \rightarrow
(\mathcal C,
\varphi)$$
 such that for all $a \in \mathcal A$ and $b \in \mathcal B$, 
$h(a \otimes b) = ab$.
\end{enumerate}
\end{proposition}

\begin{proof}  Straightforward.  Compare ~\cite{Franz what is
independence}, Proposition 3.7.
\end{proof}

\subsection{Rational rotation algebras}

A rotation algebra $\mathcal R$ is a unital algebra generated by
two invertible elements $u, v$ satisfying
\begin{equation} \label{equation rotation relation}
v u = q uv
\end{equation}
for some complex number $q$ of modulus 1.  (The $C^*$--algebra version,
which has been much studied, is a unital $C^*$--algebra generated by
two {\em unitaries} $u, v$ satisfying the relation (\ref{equation
rotation relation}); see ~\cite{Boca}.)  A rotation algebra $\mathcal R$
has a basis consisting of (certain) monomials $u^n v^m$.  We  endow
$\mathcal R$ with with the unital linear functional $\varphi$ determined
by
$$
\varphi(u^m v^n) = \begin{cases}
0 & u^m v^n \ne \one \\
1 & u^m v^n = \one. \\
\end{cases}
$$

Consider the algebra $\mathcal L = \C[[x, x\inv]]$ of formal Laurent
series, with the linear functional $\tau : \sum a_i x^i \mapsto a_0$.
Give the {\em linear space} $\L \otimes \L \cong \C[[x, x\inv, y, y\inv]]$
the algebra structure determined by
\begin{equation*}
x^m y^n x^{m'}y^{n'} = q^{m'n}x^{m + m'}y^{n + n'}.
\end{equation*}
Denote this algebra by $\L \otimes_q \L$.  The linear functional
$\tau \otimes \tau$  on $ \L \otimes_q \L$ satisfies 
$$
\tau \otimes \tau(x^m y^n) = \begin{cases}
0 & (m,n) \ne (0,0)\\
1 &  (m,n) = (0,0). \\
\end{cases}
$$
For any rotation algebra $\mathcal R$, with parameter $q$, there is a
unique unital algebra homomorphism $\psi : \L \otimes_q \L \rightarrow R$
such that $\psi(x) = u$ and $\psi(y) = v$. Moreover,
$\varphi\circ \psi = \tau\otimes \tau$; that is, $\psi$ is 
a morphism of noncommutative probability spaces.

The algebra $\L$ has a unital algebra automorphism $\gamma$ determined by
$\gamma: x^m \mapsto q^m x^m$ and satisfying 
$\tau \circ \gamma = \tau$.  The tensor square of
$\gamma$ is an algebra automorphism of $\L \otimes_q \L$; in fact, $\gamma
\otimes
\gamma$ is the inner automorphism $\gamma \otimes
\gamma = \Ad(x\inv y)$, as is easily checked.

Now suppose that $q$ is a primitive $n\th$ root of unity; in this case,
the automorphism $\gamma$ of $\L$ has order $n$, and thus
$(\L, \tau, \gamma)$ is a \ngraded \ncps.   Furthermore,
$(\L \otimes_q \L, \tau \otimes \tau, \gamma\otimes\gamma)$ is the
$\Z_n$--graded tensor square of $(\L, \tau, \gamma)$.

Consider any rotation algebra $\mathcal R$ generated by invertible
elements $u, v$ satisfying (\ref{equation rotation relation}), with
$q$ a primitive $n\th$ root of unity.  Such an algebra is called a {\em
rational rotation algebra}.  $\mathcal R$ has the inner automorphism
$\alpha = \Ad(u\inv v)$, which is of order $n$ and which satisfies
$\varphi \circ \alpha = \alpha$.  Thus, $(\mathcal R, \varphi, \alpha)$
is a  \ngraded \ncps.  

\begin{proposition} There is a unique homomorphism of \ngraded 
noncommutative probability spaces 
$\psi :  (\L \otimes_q \L, \tau \otimes \tau, \gamma\otimes\gamma)
 \rightarrow (\mathcal R, \varphi, \alpha)$  satisfying
$\psi(x) = u$ and $\psi(y) = v$.  Furthermore, the subalgebras $\A$
generated by $u$ and $\B$ generated by $v$ are \ngraded independent
in $\mathcal R$.
\end{proposition}

\begin{proof}
Straightforward.
\end{proof}

\begin{remark}  Unfortunately, the computation of the moments of the
Harper operator $u + u\inv + v + v\inv$, carried out in ~\cite{Mingo Nica}
for the universal rotation algebra with parameter $q = \pi$, cannot be
repeated here for $q$ an arbitrary root of unity, using the techniques of
this paper.  The reason is that for $n > 2$, the elements 
$u + u\inv$ and $v + v\inv$ are no longer homogeneous.
\end{remark}

\subsection{Generalized Clifford algebras}
\def\c #1 #2{C_{#1}^{(#2)}}

Let $q$ be a primitive $n\th$ root of unity.
The generalized Clifford algebra $\c m n$ (~\cite{Morris 1},
~\cite{Morris 2}) is the universal unital algebra with generators
$e_i$, ($1 \le i \le m$) subject to the relations
$$
e_i^n = \one, \quad e_j e_i = q e_i e_j \ \ \text{for}\ i < j.
$$
For $m = 2$, this is just the rotation algebra
$$
\c 2 n \cong \C[\Z_n^*] \otimes_q \C[\Z_n^*].
$$  For $m > 2$, it is the 
$m$--fold graded tensor power of $\C[\Z_n^*]$,
$$
\c m n \cong \C[\Z_n^*] \otimes_q \C[\Z_n^*] \otimes_q \cdots \otimes_q
\C[\Z_n^*].
$$
Let $\varphi$ denote  the unital linear functional $\varphi$ on 
$\c m n$ whose value on nontrivial monomials in the $e_i$ is zero.
For any $k$, the two subalgebras
 generated by  $\{e_i : i \le k\}$ and
$\{e_i: i > k\}$  are
\ngraded independent with respect to $\varphi$.

\medskip
\newsect{Crossings and ordered partitions.}

This section contains the combinatorial results on crossings of set
partitions which underlie the main results on additivity of the
graded R transform.

We recall the concepts of a {\em crossing} and of the {\em restricted
crossing number} ~\cite{Mingo Nica}   of 
a set partition of an interval $[N] = \{1, 2, \dots, N\}$.   The
restricted crossing number plays a key role in the combinatorics of set
partitions related to  ideas of independence in noncommutative
probability.

\medskip
\begin{definition} Let $P$ be a partition of $[N]$.  A {\em crossing} of
$P$ is a quadruple $a_1 < b_1 < a_2 < b_2$, where $a_i \in A$, 
$b_i \in B$, and $A$, $B$ are distinct blocks of the partition.
The number of crossings of $P$ is the {\em crossing number} 
$c(P)$.
\end{definition}

\medskip
Note that if one draws the points of $[N]$ in cyclic order on a circle
and connects each pair of points belonging to the same block by a straight
line, then  crossings corresponds to  pairs of crossing lines.    The
crossing number $c(P)$ is therefore invariant under cyclic permutations
of the points of $[N]$.  A partition $P$ is called noncrossing if
it has no crossings, $c(P) = 0$.

\medskip
\begin{definition} \label{definition reduced crossing number}
A crossing $a_1 < b_1 < a_2 < b_2$ is {\em left reduced} if
$a_1$ and $b_1$ are minimum in the blocks which contain them.
The number of left reduced crossings is the {\em restricted crossing
number} of $P$, denoted $c_0(P)$.   
\end{definition}

\medskip
Note that a partition is noncrossing if, and only if, its restricted
crossing number is zero. The restricted crossing number does not share the
symmetries of the crossing number; it is not invariant under cyclic
permutations of 
$[N]$.  For this reason, it seems a less natural statistic on set
partitions than the full crossing number.    Nevertheless, our main
technical results involve the restricted crossing number rather than the
full crossing number; these results concern the evaluation 
of quantities of
the form $\sum_{P \in \mathcal X} q^{c_0(P)}$, where $\mathcal
X$ is a collection of set partitions of an interval, and $q$ is a
parameter.  It seems that the results cannot be revised to remain valid
for the full crossing number in place of the restricted crossing number.

We will evaluate $c_0(P)$ by a sum, as follows.
For disjoint
subsets
$A, B
$ of
$[N]$, set
$$
c_0(A,B) = \begin{cases}
0 &\text{if}\quad \min B < \min A \cr
\card {\{(a,b) \in A \times B :  \min B < a < b  \}}
&\text{if}\quad
\min A <
\min B. 
\end{cases}
$$
Then for a partition $P = \{A_1, \dots, A_s\}$ of $[N]$, we have
$$c_0(P) = \sum_{i \ne j} c_0(A_i, A_j).$$

It is convenient to have
the auxiliary notions of an {\em ordered set partition} and
the {\em sorting number} of an ordered set partition.

\medskip
\begin{definition} An {\em ordered set partition} $P = (A_1, \dots, A_s)$
of a set $S$ is a sequence of subsets of $S$ such that
$\{A_1, \dots, A_s\}$ is a partition of $S$.
\end{definition}

\medskip
An ordered set partition $P$ with $s$ parts determines a surjective map
$f_P : S \rightarrow [s]$; \ namely,  $f(x) = j$ if, and only if, 
$x \in A_j$.  Conversely,  a surjective map $f: S \rightarrow [s]$
determines the ordered set partition $(f\inv(1), \dots, f\inv(s))$.

\medskip
\begin{definition} \label{definition x(P)} Given an ordered set partition
 $P = (A_1, \dots, A_s)$ of an interval $[n]$, define 
the {\em sorting number of $P$} to be $$x(P) =
\sum_{i < j}  \cardinality
{\{(x,y) \in  A_i \times A_j : x < y \}} .$$ 
Note that if $f = f_P$,
then  
$$x(P) = \card {\{(x, y) \in [n] \times [n] : x < y \text{ and } 
f(x) < f(y) \}}.
$$
\end{definition}

\medskip
\begin{lemma} \label{lemma sum 1}  Let $P = (A, B)$ be an ordered
partition of 
$[n]$ with two parts.  Let $q$ be a proper $n\th$ root of unity.
Let $\sigma$ denote the cyclic permutation of $[n]$,
$\sigma = (1, 2, \dots, n)$, in cycle notation.
Then
$$
\sum_{0 \le k \le n-1} q^{x(\sigma^k(P))} = 0.
$$
\end{lemma}

\begin{proof}  Let $r = \card A $.  \ If $n \in B$, then
$x(\sigma(P)) = x(P) - r$.  If $n \in A$, then
$x(\sigma(P)) = x(P) + (n-r) \equiv x(P) -r  \mod n$.
It follows that for all $k$, $x(\sigma^k(P)) \equiv x(P) - kr \mod n$.
Therefore,
$$
\sum_{0 \le k \le n-1} q^{x(\sigma^k(P))} = q^{x(P)} \sum_{0 \le k \le
n-1} q^{- k r} = 0.
$$
\end{proof}
\medskip

\begin{definition} \label{definition cyclic action 1}
We define an action of $\Z_n$ on ordered set partitions of $[n]$, as 
follows
Let $P = (A_1, \dots, A_s)$.  Put $B = [n] \setminus A_1$.
Put $\sigma(j) = j+1  \mod n$.  Put $B' = \sigma(B)$, and let
$\theta : B \rightarrow B'$ be the unique order preserving bijection.
Define
$
\mathbold \sigma(P) = (A'_1, \dots, A'_s),
$
where $A'_1 = \sigma(A_1)$, and $A'_j = \theta(A_j)$, for $j > 1$.
\end{definition}

\begin{remark} It is easy to see that $\mathbold \sigma$ has order $n$
and that $P$ is fixed if, and only if, $P$ has only one part.
\end{remark}

\begin{example} $n = 7$, $P = (A_1, A_2, A_3)$, where
$A_1 = \{1,5,6\}$, $A_2 = \{2, 7  \}$, and $A_3 = \{3,4  \}$.
Then $\mathbold \sigma(P) = \{A'_1, A'_2, A'_3\}$, where
$A'_1 = \{ 2, 6,7 \}$, $A'_2 = \{1,5   \}$, $A'_3 = \{ 3, 4  \}$.
The map $f : [n] \rightarrow [3]$ corresponding to $P$ is
$\begin{pmatrix} 1 & 2 & 3 & 4 & 5 & 6 & 7 \cr 1 & 2 & 3 & 3 & 1 & 1 & 2
\end{pmatrix}$.  The map corresponding to $\mathbold \sigma(P)$ is
$\begin{pmatrix} 1 & 2 & 3 & 4 & 5 & 6 & 7 \cr 2 & 1 & 3 & 3 & 2 & 1 & 1 
\end{pmatrix}$.
\end{example}

\medskip
\begin{lemma} \label{lemma sum 2}  Let $P = (A_1, \dots, 
A_s)$ be an ordered partition of 
$[n]$ with $s \ge 2$ parts.  Let $q$ be a proper $n\th$ root of
unity. Let $\mathbold \sigma$ denote the generator of the $\Z_n$ action
on ordered partitions of $[n]$.
Then
$$
\sum_{0 \le k \le n-1} q^{x(\mathbold \sigma^k(P))} = 0.
$$ 
\end{lemma}

\begin{proof} For two subsets $A, B$ of $[n]$, let $x(A, B) =
\card {\{(x, y) \in A \times B : x < y  \}}$ .
Write $\mathbold \sigma(P) = (A'_1, \dots, A'_s)$.
Set $B = [n] \setminus A_1$,  and $B' = [n] \setminus A'_1 =
\sigma(B)$.  We have
$$
x(P) = \sum_{i < j} x(A_i, A_j) = x(A_1, B) + \sum_{2 \le i < j} x(A_i,
A_j).
$$
Similarly,
$$
x(\mathbold \sigma(P)) = \sum_{i < j} x(A'_i, A'_j) = x(A'_1, B') +
\sum_{2
\le i < j} x(A'_i, A'_j).
$$
But $x(A'_i, A'_j) = x(A_i, A_j)$ if $2 \le i < j$.
Note that $(A_1, B)$ is an ordered partition with 2 parts and
$\sigma(A_1, B) = \mathbold \sigma(A_1, B) = (A'_1, B')$.  By the
argument of Lemma \ref{lemma sum 1}, 
$x(A'_1, B') \equiv x(A_1, B) - r  \mod n$, where $r = \card {A_1} $.
Therefore $q^{x(\mathbold \sigma(P))} = q^{x(P)} q^{-r}$.  The conclusion
follows.
\end{proof}

\medskip
The family of all set partitions of a set $S$ will be denoted by 
$\P(S)$.    Given sets $T \subseteq S$ and a set partition $P$ of
$S$,  the partition $P\rest T$ induced on $T$ by $P$ is
$\{ T \cap A : A \in P \text{ and } T \cap A \ne
\emptyset \}$.
   
Fix natural numbers $n$ and $m$.  We consider set partitions
of $[m n]$ such that each part has cardinality divisible by $n$.
Denote the collection of such set partitions by $\P_n([mn])$.
For $0 \le k \le m-1$, let $J_k$ denote the interval
$\{k n +1, \dots, (k+1) n\}$.  Let $\P_n^0([mn])$ denote
those set partitions $P \in \P_n([mn])$ such that
for all $k$ and for all $A \in P$, either $J_k \subseteq A$, or
$J_k \cap A = \emptyset$.  Equivalently, the requirement on
$P$ is that for each $k$, the set partition
$P\rest {J_k}$  induced on $J_k$ by $P$ is the trivial set partition $\{
J_k
\}$.

We are going to define an action of $\Z_n$ on $\P_n([mn])$
whose fixed point set is precisely $\P_n^0([mn])$.  
If $P \in  \P_n^0([mn])$, set $\mathbold \sigma(P) = P$.

For $P \in (\P_n([mn])  \setminus \P^0_n([mn]))$, there
exists
$k$ ($1 \le k \le m$) such that the induced  partition  $P\rest{J_k}$  is
not the trivial partition of
$J_k$.  Let $k_0$ be the maximum of such $k$.

Write $P = (A_1, \dots, A_\ell)$, with the order given by
$i < j$ \ if, and only if, $\min A_i < \min A_j$.
Consider the induced partition on $J_{k_0}$, with the induced order; 
namely,   $P\rest{J_{k_0}} = (B_1, \dots, B_s)$, with
$B_j = A_{i_j} \cap J_{k_0}$, and $i_1 < i_2 < \dots < i_s$.
By definition of $k_0$, we have $s \ge 2$.

We let  $\Z_n$ act on $P$ via its action on
$P\rest {J_{k_0}}$.  More precisely, let the generator $\mathbold\sigma$
act on $P\rest {J_{k_0}}$ as in Definition \ref{definition cyclic action
1}, and  put
$\mathbold \sigma(P\rest {J_{k_0}}) = (B'_1, \dots, B'_s)$.
Define $$\mathbold \sigma(P) = (A'_1, \dots, A'_\ell),$$ 
where
$$
A'_i = \begin{cases}
A_i &\text{if}\quad A_i \cap J_{k_0} = \emptyset, \cr
 B'_j \cup  \bigcup_{k \ne k_0} (A_i \cap J_k) &\text{if}\quad
 A_i \cap J_{k_0} = B_j.
\end{cases}
$$

Evidently, $\mathbold \sigma$  has order $n$, and $P$ is fixed by 
$\mathbold \sigma$ if, and only if, $P \in \P_n^0([mn])$.

\begin{proposition} \label{proposition sum 3}
Let $P \in (\P_n([mn]) 
\setminus
\P^0_n([mn]))$, and let $q$ be a nontrivial $n\th$ root of unity.
Then
$$
\sum_{0 \le k \le n-1} q^{c_0(\mathbold \sigma^k(P))} = 0.
$$ 
\end{proposition}

\begin{proof} Write $P = (A_1, \dots, A_t)$, with
$i < j$ if, and only if, $\min A_i < \min A_j$.  Write
$\mathbold \sigma(P) = (A'_1, \dots, A'_t)$.  Let
$k_0$ be the maximum of those $k$ such that the restricted
partition $P\rest{J_k} $ is nontrivial.

Set $c_0(A_i, A_j)_{k,\ell} =$
$$
\card {\{(a,b) : a \in A_i \cap J_k, \ b \in A_j \cap J_\ell,\text{ and } 
\min A_i < \min A_j < a < b   \}
} .$$
Then $$c_0(A_i, A_j) = \sum_{k \le \ell} c_0(A_i, A_j)_{k , \ell}.$$
Observe that
$c_0(A_i, A_j)_{k , \ell} = c_0(A'_i, A'_j)_{k , \ell}$, unless
$k = \ell = k_0$ and  both $A_i \cap J_{k_0}$ and \break
$A_j \cap J_{k_0}$ are nonempty.

It follow that $$ c_0(\mathbold \sigma(P)) - c_0(P)   =
\sum_{i < j} c_0(A'_i, A'_j)_{k_0, k_0} - 
\sum_{i < j} c_0(A_i, A_j)_{k_0, k_0}
.$$

Note that if $A \in P$ and $A \cap J_{k_0} \ne \emptyset$, then
$\min A < k_0 n + 1$.  This is because 
$\card A \equiv 0  \mod n$ by definition of $\P_n$, 
 $\card {A \cap J_k} \equiv 0  \mod n$ if $k > k_0$,  by definition of
$k_0$, and
$\card {A \cap J_{k_0}} \not\equiv 0  \mod n$ by assumption.  This implies
there exists $k < k_0$ such that $A \cap J_k \ne \emptyset$.
It follows from this observation that if $A_i \cap J_{k_0}$ and
$A_j \cap J_{k_0}$ are nonempty, then
$$
c_0(A_i, A_j)_{k_0, k_0} = x(A_i \cap J_{k_0}, A_i \cap J_{k_0}),
$$
and, therefore,
$$
\sum_{i < j} c_0(A_i, A_j)_{k_0, k_0} = \sum_{i < j} x(A_i \cap J_{k_0},
A_i \cap J_{k_0}) = x(P\rest {J_{k_0}}),
$$
and similarly for $\mathbold \sigma(P)$ replacing $P$.  It follows  that
$$
c_0(\mathbold \sigma(P)) - c_0(P) =
x(\mathbold \sigma(P)\rest{J_{k_0}}) -
x(P\rest{J_{k_0}}),
$$
and the conclusion of the proposition follows from this and Lemma
\ref{lemma sum 2}.
\end{proof}

\medskip
\begin{corollary} \label{corollary q sum 4.5} Fix natural numbers $m$ and
$n$. For $q$  a proper
$n\th$ root of unity,
$$
\sum_{P \in \P_n[mn]} q^{c_0(P)} = \card {\P^0_n[mn]} .
$$
\end{corollary}

\begin{proof}  It is easy to see that for $P \in {\P^0_n[mn]}$,
$c_0(P) \equiv 0  \mod n$, so $q^{c_0(P)} = 1$.  On the other hand,
if $P \in (\P_n[mn] \setminus \P^0_n[mn])$, then
$P$ belongs to a nontrivial orbit $\mathcal O$ for the $\Z_n$ action, and,
by Proposition \ref{proposition sum 3}, the sum over the orbit is zero,
$\sum_{Q\in \mathcal O}q^{c_0(Q)} = 0$.
\end{proof}

\medskip
\begin{definition}  Let $\la$ be a {\em partition of the natural number}
$N$; that is, $\la$ is a sequence of natural numbers
$\la_1 \ge \la_2 \ge \cdots \ge \la_r $ such that $\sum_i \la_i = N$.
The $\la_i$ are called the {\em parts} of $\la$.
A set partition $P$ of $[N]$ is said to be {\em of type $\la$} if
the cardinalities of the parts of $P$, listed in decreasing order,
are the parts of $\la$.  If $\la$ is a partition of $N$ and $m$ is a
natural number,  the partition $(m \la_1, m\la_2, \dots, m \la_r)$ 
of $m N$ is denoted by $m \la$.
\end{definition}

\medskip
\begin{corollary} \label{corollary q sum 5} Fix natural numbers $m$ and
$n$. Let $\la$ be a partition of $m$.
For $q$  a proper
$n\th$ root of unity,
$$
\sum_{\begin{array}{c} 
{\scriptstyle P \in \P_n[mn]} \cr
{\scriptstyle {\rm type}(P) = n \la}
\end{array}
}  
q^{c_0(P)} \  = \  \card {\{ P \in \P[m]  : {\rm type}(P) = \la
\}} .
$$ 
\end{corollary}

\begin{proof}  The set of $P \in \P_n[mn]$ that are of type
$n \la$  is invariant under the $\Z_n$ action on  $\P_n[mn]$.
The proof of the previous corollary shows that  
$$\sum_{\begin{array}{c} 
{\scriptstyle P \in \P_n[mn]} \cr
{\scriptstyle {\rm type}(P) = n \la}
\end{array}
}  
q^{c_0(P)}\  = \ 
 \card {\{ P \in \P_n^0[mn]  : {\rm type}(P) = n\la
\}} .
$$
But ${\{ P \in \P_n^0[mn]  : {\rm type}(P) = n\la
\}}   $  is in one to one correspondence with \break 
${\{ P \in \P[m]  : {\rm type}(P) = \la
\}}$.
\end{proof}

\medskip
\medskip
\begin{remark}  For $n = 2$ and $P \in \P_n[nm]$, Mingo and Nica show that
$c_0(P) \cong c(P) \mod n$; therefore for $n = 2$, Corollaries
\ref{corollary q sum 4.5} and \ref{corollary q sum 5} are  also
valid with the restricted crossing number replaced by the full crossing
number.  However, for $n \ge 3$, the congruence $c_0(P) \cong c(P) \mod
n$ does not hold, and the analogues of Corollaries
\ref{corollary q sum 4.5} and \ref{corollary q sum 5} with $c_0(P)$
replaced by $c(P)$ are not valid, as shown by computations.
\end{remark}
\medskip

\newsect{$\Z_n$-graded noncommutative probability spaces.}

Let $(\A, \varphi, \gamma)$ be a \ngraded \ncps.
If $a$ is homogeneous, say $a \in A_q$,
then the only moments  $\varphi(a^k)$ which are (possibly) nonzero are
those for which
$a^k \in A_1$, or, equivalently, $q^k =1$.
This is true if, and only if, $k$ divides the order of $q$ in
$\Z_n^*$, which is also the order of $\gamma
\rest {A_q}$.

\medskip
\begin{lemma} Let $n$ be a natural number.  Suppose $\mu$ satisfies
$\mu(X^k) = 0$ unless $n$ divides $k$.   For any $q$, the sequence
$(\a q k)$ of $q$-cumulants of $\mu$ also satisfies $\a q k = 0$
unless $n$ divides $k$.
\end{lemma}

\begin{proof}  This is easily seen by induction on $k$.
\end{proof}

\medskip
\begin{definition} Let $n$ be a natural number and let $q$ be a proper 
$n\th$ root of unity.  Let $\mu : \C[X] \rightarrow \C$ be a linear
functional with $\mu(1) = 1$ and $\mu(X^k) = 0$ unless $n$ divides $k$.
Let $(\a q k)_{k \ge 1}$ be the sequence of $q$--cumulants of $\mu$.
The {\em $\r n q $-- transform of $\mu$} is the power series
$$
\r n q[\mu](z) = \sum_{k \ge 1} \frac{\a q {nk} }{k !} \ z^{nk}.
$$
For an arbitrary linear functional $\mu : \C[X] \rightarrow \C$ with
$\mu(1) = 1$, let $(\a 1 k )_{k \ge 1}$ be the sequence of $1$--cumulants
of
$\mu$. The $R_1$--transform of $\mu$ is the power series
$$
R_1[\mu](z) = \sum_{k \ge 1} \frac{\a 1 k }{k !} \ z^k.
$$

\end{definition}

\medskip
\begin{remark} It is necessary to specify $n$ as well as $q$ in 
$\r n q$, as it is
not assumed that $q$ is a {\em primitive} $n\th$ root of unity.
\end{remark}

\medskip

\begin{theorem}  \label{theorem alpha = beta}  Let $n$
be a natural number and let  $q$ be a proper $n\th$ root of unity.
 Let
$(\A,
\varphi)$ be an algebraic noncommutative probability space. Suppose
$a \in \A$ satisfies \break $\varphi(a^k) = 0$ unless $n$ divides $k$.
Let $(\a q k)_{k \ge 1}$ be the sequence of $q$--cumulants of
$\mu_a$.  (Recall that  $\a q k = 0$ unless $n$ divides $k$.)  Let
$(\beta_k)_{k \ge 1}$ be the sequence of $1$--cumulants of $\mu_{a^n}$.
Then for all $k \ge 1$, $$\a q {nk} = \beta_k.$$
\end{theorem}

\begin{proof}  The proof is the same, {\em mutatis mutandis}, as that of
Corollary 1.8 in ~\cite{Mingo Nica}.  We give the proof here for the sake
of completeness.  The proof goes by induction on $k$.  For $k = 1$, we
have  $$\a q  {n} = \mu_a(X^n) = \varphi(a^n) = \mu_{a^n}(X) =
\beta_1.$$
Fix $k > 1$ and suppose the equality $\a q  {n k'} = \beta_{k'}$
holds for all $k' < k$.  We have
$$
\varphi(a^{kn}) = \mu_a(X^{kn}) = \sum_{P \in \P[kn]} q^{c_0(P)}
\prod_{B \text{ a block of } P} \a q  {\card {B} }.
$$
Since $\a q  r = 0$ unless $n$ divides $r$, the latter sum reduces to
$$
\sum_{P \in \P_n[kn]} q^{c_0(P)}
\prod_{B \text{ a block of } P} \a q  {\card {B} } =
\sum_{\la \text{ partition of } k}
 (\sum_{ \begin{array}{c}
{\scriptstyle P \in \P_n[kn]} \cr
{\scriptstyle P \text{ of type } n \la}
\end{array}}  
q^{c_0(P)} )
\prod_{i} \a q  {n \la_i }.
$$
Applying Corollary \ref{corollary q sum 5} to the inner sum, we obtain
\begin{equation} \label{equation alpha = beta 1}
\varphi(a^{kn}) = 
\sum_{\la \text{ partition of } k}
 (
\sum_{ \begin{array}{c}
{\scriptstyle P \in \P[k]} \cr
{\scriptstyle P \text{ of type }  \la}
\end{array}
}  
1)\prod_{i} \a q  {n \la_i }.
\end{equation} 
On the other hand, we have
\begin{equation} \label{equation alpha = beta 2}
\begin{aligned}
\varphi(a^{kn}) &= \mu_{a^n}(X^{k}) =
\sum_{P \in \P[k]} 
\prod_{B \text{ a block of } P} \beta_{\card {B} } \cr
&= \sum_{\la \text{ partition of } k}
 (
\sum_{ \begin{array}{c}
{\scriptstyle P \in \P[k]} \cr
{\scriptstyle P \text{ of type }  \la}
\end{array}
}  
1)\prod_{i} \beta_{ \la_i }.
\end{aligned}
\end{equation}
Comparing Equations (\ref{equation alpha = beta 1}) and
(\ref{equation alpha = beta 2}) gives
\begin{equation}
\sum_{\la \text{ partition of } k}
 (
\sum_{ \begin{array}{c}
{\scriptstyle P \in \P[k]} \cr
{\scriptstyle P \text{ of type }  \la}
\end{array}
}  
1) \left[\prod_{i} \beta_{ \la_i } - \prod_{i} \a q  {n \la_i} \right] =
0.
\end{equation}
All of the expressions $\left[\prod_{i} \beta_{ \la_i } - \prod_{i}
\a q  {n \la_i} \right]   $ are zero, by the induction assumption,
with the exception of that corresponding to $\la =  (k)$,  which is
$\beta_k - \a q  {nk}$.  It follows that 
$\beta_k - \a q  {nk} = 0$ as well.
\end{proof}

\medskip
\begin{corollary} \label{rq transform and r1 transform} Let $n$
be a natural number and let  $q$ be a proper $n\th$ root of unity.
 Let
$(\A,
\varphi)$ be an algebraic noncommutative probability space. Suppose
$a \in \A$ satisfies  $\varphi(a^k) = 0$ unless $n$ divides $k$.
Then
$$
\r n q[\mu_a](z) = R_1[\mu_{a^n}](z^n).
$$
\end{corollary}

\begin{proof}  This is immediate from the theorem and the definition of
the $R$--transforms.
\end{proof}

\newsect{$\Z_n$--graded independence and the $\Z_n$--graded R--transform.}

Let $n$ be a natural number and let $(\mathcal C, \varphi, \gamma)$ be an
 $\Z_n$--graded algebraic noncommutative probability space.

\medskip
\begin{lemma} \label{lemma freshman power rule}   Let $\A$ and
$\mathcal B$  be {\em$\Z_n$--graded independent} subalgebras of
$\mathcal C$. Let $q$ be a primitive $n\th$ root of unity such
that $\A$ and \ 
$\mathcal B$ $q$--commute.
 Let $a \in
\A$ and
$b
\in
\mathcal B$ be homogenous elements of the same $q$--degree $r$.  Write $n'
= n/{\gcd(r,n)}$. Then  $(a + b)^{n'} = a^{n'} +b^{n'}$.
\end{lemma}

\begin{proof} The elements $a$ and $b$ satisfy
$b a =
q^{r^2} a b$.  Write $q'$ for $q^{r^2}$, and  note that $q'$ is
a primitive ${n'}\th$ root of unity.  Let $0 < k < n'$, and consider the 
sum $S$ of 
those terms in the expansion of $(a + b)^{n'} $ that are of degree
$k$ in $a$ and degree $n' - k$ in $b$.  It follows from the relation
$b a =
q' a b$ that 
$$
S = a^k b^{n'-k} \sum_{(A,B)} ({q'})^{x((B, A))},
$$
where the sum is over all ordered set partitions $(B, A)$ of $[n']$
with two parts of cardinalities $\card {B} = n'- k$, $\card {A} = k$,
and ${x((B, A))} $  is as in Definition \ref {definition x(P)}.
Consider the $\Z_{n'}$ action on such ordered partitions, induced by
the action on $[n']$ by powers of the cyclic permutation
$(1, 2, \dots, n')$.  According to Lemma \ref{lemma sum 1}, for any
orbit $\mathcal O$,  the sum over the orbit
$
\sum_{P \in \mathcal O} (q')^{x(P)}
$
is equal to zero, since $q'$ is a proper ${n'}\th$ root of unity.  But
then
$$
\sum_{(A,B)} ({q'})^{x((B, A))} = \sum_{\mathcal O} \sum_{P \in \mathcal
O} (q')^{x(P)} = 0.
$$
\end{proof}

\medskip
\begin{definition} \label{definition graded r transform}
 Let $q$ be a fixed primitive $n\th$ root of unity.   Let $a \in \A$ be
homogeneous.   If $a \not\in \A_1$, put $n' = n/{\gcd(\delta_q(a),n)}$.
 The {\em
graded $\r q n$--transform}  $a $ is
defined by
$$
\r n q[a](z) = 
\begin{cases}
R_1[\mu_a](z) & \text{if } a \in \A_1 \\
\r  {n'}  {q^{r}}[\mu_a](z) & \text{if } a \in \A_{q^r} \ne \A_1.\\
\end{cases}
$$
\end{definition}
\medskip

The proofs of the following two results follow those of the corresponding
results in ~\cite{Mingo Nica}.

\medskip
\begin{theorem} \label{theorem formula for rq transform} Let $a$ be a
homogeneous element of
$\A$, and let $q$ be a primitive $n\th$ root of unity.   Set 
$$n' =
\begin{cases}
1 & a \in \A_1 \\
n/{\gcd(\delta_q(a),n)} & \text{otherwise}.
\end{cases}
$$  The graded
$\r n q$--transform of
$a$ satisfies
$$
\r n q[a](z) = \log\left[ \sum_{k \ge 0} \frac{\mu_a(X^{k n'})}
{k !} \  z^{k n'}  
\right].
$$
\end{theorem}

\begin{proof} If $a \in A_1$, this is immediate from Equation (\ref{equation formula
classical R transform}).
 Otherwise, put $r = \delta_q(a)$.   We have
$$\begin{aligned}
\r n q[a](z) &= \r  {n'}  {q^{r}}[\mu_a](z) = 
R_1[\mu_{a^{n'}}](z^{n'}) \cr
&= \log\left[ \sum_{k \ge 0} \frac{\mu_{a^{n'}}(X^k)}{k !} \ (z^{n'})^k
\right ] \cr
&= \log\left[ \sum_{k \ge 0} \frac{\mu_{a}(X^{k n'})}{k !} \
z^{kn'}
\right ],
\end{aligned}
$$
where the second equality follows from Corollary \ref{rq transform and
r1 transform} and the third from Equation (\ref{equation formula
classical R transform}).
\end{proof}

\medskip

\begin{theorem} \label{theorem rq transform linearizes} Let $\A$ and
$\mathcal B$  be {\em$\Z_n$--graded independent} subalgebras of $\mathcal
C$. Let $a \in \A$ and $b \in \mathcal B$ be homogenous
elements of the same  degree.  Then
$$
\r  n q [a + b] = \r n q [a]  + \r n q [b].  
$$
\end{theorem}

\begin{proof} Note that the fixed point algebras $\A_1$ and $\mathcal
B_1$ are classically independent; that is, they mutually commute, and $\varphi(x y) = \varphi(x)\varphi(y)$
for all $x \in \A_1$, $y \in \mathcal B_1$.  Furthermore,
for
classically independent random variables $x$ and $y$, one has
\begin{equation} \label{equation classical r transform linearizes}
R_1[\mu_{x + y}] = R_1[\mu_x] + R_1[\mu_y].
\end{equation}
Therefore,  for $a, b \in \A_1$, the result follows from 
Definition \ref{definition graded r transform}  and Equation
(\ref{equation classical r transform linearizes}).

If $a, b \not\in \A_1$, let $r$ denote their $q$-degree and put  $n' =
n/{\gcd(r,n)}$.  We have
$$
\r n q [a + b](z) = \r {n'} {q^r}[\mu_{a + b}](z) = R_1[\mu_{(a
+ b)^{n'}}](z^{n'}),
$$
where the first equality is by Definition  \ref{definition graded r
transform} and the second by Corollary \ref{rq transform and r1
transform}.

According to Lemma \ref{lemma freshman power rule}, we have
$(a +b)^{n' }  = a^{n' } + b^{n' }  $, so 
$$
 R_1[\mu_{(a
+ b)^{n'}}](z^{n'})= R_1[\mu_{ a^{n' } + b^{n' }}](z^{n'}).
$$
Observe that $a^{n'} $ and $b^{n'}$  are elements of  the classically
independent fixed point algebras
$\A_1$, $\mathcal B_1$, so
$$
\begin{aligned}
R_1[\mu_{ a^{n' } + b^{n' }}](z^{n'}) &=
R_1[\mu_{ a^{n' } }](z^{n'})  +
R_1[\mu_{ b^{n' }}](z^{n'})  \cr
&=\r {n'} {q^r}[\mu_{a }](z) +
\r {n'} {q^r}[\mu_b](z) \cr
&= \r n q [a ](z)   + \r n q [ b](z),
\end{aligned}
$$
using Equation (\ref{equation classical r transform linearizes}),
Corollary \ref{rq transform and r1 transform}   and Definition
\ref{definition graded r
transform}.
\end{proof}

\end{document}